\documentclass[pdflatex,sn-mathphys-num]{sn-jnl}

\usepackage{graphicx}%
\usepackage{multirow}%
\usepackage{amsmath,amssymb,amsfonts}%
\usepackage{amsthm}%
\usepackage{mathrsfs}%
\usepackage[title]{appendix}%
\usepackage{xcolor}%
\usepackage{textcomp}%
\usepackage{manyfoot}%
\usepackage{booktabs}%
\usepackage{algorithm}%
\usepackage{algorithmicx}%
\usepackage{algpseudocode}%
\usepackage{listings}%
\usepackage{subfigure}
\newcommand\onenorm[1]{\left\lvert#1\right\rvert}
\newcommand\twonorm[1]{\left\lVert#1\right\rVert}

\begin{document}


\title[Article Title]{Multivariable Stochastic Newton-Based Extremum Seeking with Delays}


\author[1]{\fnm{Paulo Cesar} \sur{Souza Silva}}\email{cesar.paulo151@hotmail.com}
\equalcont{These authors contributed equally to this work.}

\author[1]{\fnm{Paulo César} \sur{Pellanda}}\email{pellanda@ime.eb.br}
\equalcont{These authors contributed equally to this work.}

\author*[2]{\fnm{Tiago Roux} \sur{Oliveira}}\email{tiagoroux@uerj.br}
\equalcont{These authors contributed equally to this work.}

\affil[1]{\orgdiv{Department of Defense Engineering}, \orgname{Military Institute of Engineering}, \orgaddress{\street{Praça General Tibúrcio 80, Urca}, \city{Rio de Janeiro}, \postcode{22290-270}, \state{RJ}, \country{Brazil}}}
\affil[2]{\orgdiv{Department of Electronics and Telecommunication Engineering},\\ \orgname{State University of Rio de Janeiro}, \orgaddress{\street{São Francisco Xavier 524, Maracanã}, \city{Rio de Janeiro}, \postcode{22550-900}, \state{RJ}, \country{Brazil}}}




\abstract{This paper presents a Newton-based stochastic extremum-seeking control method for real-time optimization in multi-input systems with distinct input delays. It combines predictor-based feedback and Hessian inverse estimation via stochastic perturbations to enable delay compensation with user-defined convergence rates. The method ensures exponential stability and convergence near the unknown extremum, even under long delays. It extends to multi-input, single-output systems with cross-coupled channels. Stability is analyzed using backstepping and infinite-dimensional averaging. Numerical simulations demonstrate its effectiveness in handling time-delayed channels, showcasing both the challenges and benefits of real-time optimization in distributed parameter settings.

\textcolor{white}{
This paper presents a Newton-based stochastic extremum-seeking control method for real-time optimization in multi-input systems with distinct input delays. It combines predictor-based feedback and Hessian inverse estimation via stochastic perturbations to enable delay compensation with user-defined convergence rates. The method ensures exponential stability and convergence near the unknown extremum, even under long delays. It extends to multi-input, single-output systems with cross-coupled channels. Stability is analyzed using backstepping and infinite-dimensional averaging. Numerical simulations demonstrate its effectiveness in handling time-delayed channels, showcasing both the challenges and benefits of real-time optimization in distributed parameter settings.}

}

\maketitle

\section{Introduction}\label{sec1}
Extremum seeking control (ESC) is a model-free, adaptive, and real-time optimization method \citep{c19}. It aims to determine the extremum --- either a minimum or maximum --- of a nonlinear map \citep{c1}. 
ESC operates without requiring explicit knowledge of the plant or the function to be optimized,
relying only on the existence of an extremum in the nonlinear function \citep{c2,Aminde2013}.

\textcolor{black}{Despite extensive research and numerous recent studies on extremum seeking --- spanning both theoretical advancements and practical applications \citep{c2,c3,c4,c5,c19,c20,CBA4B,ArtElly} --- the problem of multivariable Newton-based stochastic extremum seeking control in the presence of distinct input delays remains unaddressed \citep{c7,c6}. This gap is particularly relevant in applications where significant delays arise due to post-processing of the plant's measured output, leading to a substantial lag in generating the control input. This is the case in laser-based light sources for photolithography in semiconductor manufacturing, where image processing delays are inherent \cite{ref_15_TAC,ref_16_TAC}, as well as in various chemical and biochemical processes where sample analysis induces unavoidable latency. A notable example is the phase lag observed in batch cultures within bioreactors \cite{ref_17_TAC,ref_18_TAC}, which exemplifies the biological optimization process being hindered by delays. These delays are typically known, constant, distinct, and relatively large, making their compensation a critical challenge in real-time control systems.}

Furthermore, the stochastic extremum controller offers significant advantages over its deterministic counterpart, including the ability to avoid local extrema and achieve a faster convergence rate \citep{c6}. Additionally, the Newton algorithm provides advantages over the gradient algorithm \citep{c4}, such as removing the convergence rate's dependence on the unknown second derivative (Hessian) of the nonlinear map, allowing this rate to be arbitrarily assigned by the designer or user \citep{c2,c4}.

\textcolor{black}{On the other hand, while ESC generally ensures a reasonable convergence rate, the presence of delays in the closed-loop system --- if not properly compensated --- can significantly degrade performance or even lead to instability, as ESC is not inherently robust to arbitrarily long delays \citep{c7}. A major milestone in this context was presented in \cite{c7}, where the deterministic ESC framework was extended to partial differential equations (PDEs), addressing the design and analysis of multiparameter static maps subject to such delays. The study showed that these delays can be effectively modeled as first-order hyperbolic transport PDEs, laying the groundwork for further extensions to other classes of PDEs, as explored in \cite{Oliveira2022}.
}

\citet{LIUFinal} laid the foundation for multivariable Newton-based stochastic extremum-seeking control without delays, establishing local stability through convergence in probability. \textcolor{black}{Later, \citet{c7,CBA2016B,c23} extended deterministic ESC to account for input and/or output delays using gradient- and Newton-based algorithms. Building on these developments, \citet{NewHigh} introduced a stochastic ESC framework based on the Newton algorithm, specifically targeting the optimization of higher-order derivatives of unknown maps subject to output delays. Their use of stochastic perturbations allowed for the compensation of arbitrarily long output delays while optimizing a dynamic, single-input map.}

In this paper, we present a solution to the \textcolor{black}{open problem} of multivariable stochastic ESC with multiple and distinct input delays by employing predictor feedback to estimate the inverse of the unknown Hessian based on stochastic sinusoidal perturbations \citep{c6}. The stability analysis is rigorously developed using backstepping transformations \citep{c1} and infinite-dimensional averaging theory \citep{AHale,ALehman}, effectively addressing the infinite-dimensional state induced by the delays. 
This work builds on preliminary results presented in a previous conference version \cite{CODIT2024}, with \textcolor{black}{the main advancement being the development of a stochastic Newton-based ESC algorithm capable of compensating for distinct input delays. This contrasts with the previously considered stochastic gradient-based ESC method, which was limited to handling only uniform delays across all input channels.}

\bigskip
\textbf{Notation:} The partial derivatives of $u(x,t)$ are denoted by $u_t(x,t)$ and $u_x(x,t)$ or, occasionally, by $\partial_tu_{av}(x,t)$ and $\partial_xu_{av}(x,t)$ to refer to the average signal operator,  $u_{av}(x,t)$. The 2-norm of a finite-dimensional state vector $\vartheta(t)$ is denoted by single bars, $\onenorm{\vartheta(t)}$, while norms of functions (in terms of $x$) are indicated by double bars. We denote the spatial $\mathcal{L}_{2}[0,1]$ norm of the PDE state $u(x, t)$ as $\twonorm{u(t)}^{2}_{\mathcal{L}_{2}([0,1])} := \int_{0}^{1}u^2(x,t)dx$. For simplicity, we drop the subscript  $\mathcal{L}_{2}([0,1])$ in the text, so $\twonorm{\cdot} = \twonorm{\cdot}_{\mathcal{L}_{2}([0,1])}$, 
  unless otherwise specified. Consider a generic nonlinear system $\dot{x}=f(t,x,\epsilon)$, where $x \in$ $\mathbb{R}^n$, and  $f(t,x,\epsilon)$ is a periodic function with period $T$, i.e., $f(t+T,x,\epsilon)=f(t,x,\epsilon)$. For sufficiently small $\epsilon > 0$, the average model can be expressed as $\dot{x}_{av}=f_{av}(x_{av})$, with $f_{av}=\frac{1}{T}\int_{0}^{T} f(\tau,x_{av},0) d\tau$, where $x_{av}(t)$ denotes the average version of the state $x(t)$ \cite{c13}. As defined in \citep{c13}, a vector-valued function $f(t,\epsilon) \in \mathbb{R}^n$ is said to be of order $O(\epsilon)$,  within the interval $[t_1,t_2]$, if there exist positive constants $k$ and $\epsilon^{\ast}$, such that $|f(t,\epsilon)| \leq k \epsilon$, for all $\epsilon \in [0,\epsilon^\ast]$ and $t \in [t_1, t_2]$. Occasionally, 
$O(\epsilon)$ is used to denote an order of magnitude relation for sufficiently small $\epsilon$. The acronym ODEs refers to ordinary differential equations, while the acronym PDEs designates partial differential equations. 
\textcolor{black}{Given an $\mathbb{R}^n$-valued signal $f(t)$, the notation $f^D(t)$ is defined as $f^D(t):=\big[f_1(t-D_1)\; f_2(t-D_2)\; ...\; f_n(t-D_n)\big]^T$. With slight abuse of notation, we also express this operator as $f^D(t):=f(t-D)=e^{-Ds}[f(t)]$, where $D=\textrm{diag}\{D_1, D_2, \cdot \cdot \cdot, D_n \}$ represents a diagonal matrix of delays. Each input channel experiences a distinct, known, and constant delay $D_i \geq 0$, $\forall i \in \{1,2,...,n \}$, with $s$ denoting the Laplace operator.}

\section{Multivariable Stochastic Extremum Seeking with Distinct Delays}\label{sec2}
Multivariable ESC is used in scenarios where the primary objective is to maximize (or minimize) the output $y \in \mathbb{R}$ of an unknown nonlinear static map $y = Q(\theta)$ by adjusting the input vector $\theta = [\theta_1 \; \theta_2 \; ... \; \theta_n]$ in real-time.

\begin{figure}[htb!]
\centering \includegraphics[scale=0.7]{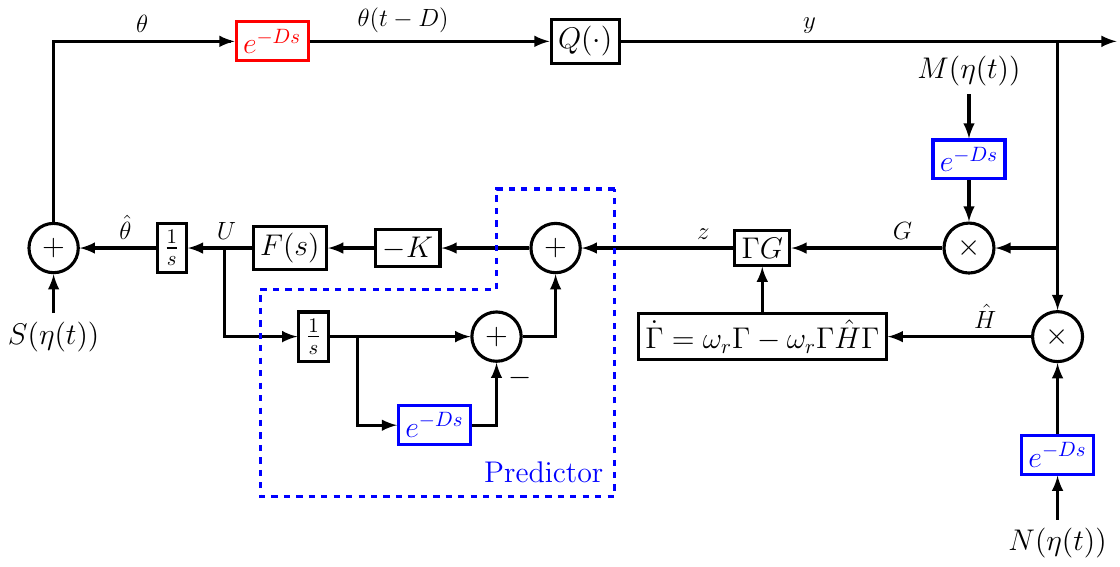}
\caption{Block diagram of the basic prediction scheme for compensating distinct input delays in stochastic extremum seeking using the Newton algorithm. The predictor feedback (\ref{eq:29}) is shown in its vector form with $F(s)=diag\{c_1/(s+c_1),...,c_n/(s+c_n)\}$ and $K=diag\{k_1,...,k_n\}$. The red block indicates the introduced delays, while the blue blocks show modifications to the classical stochastic Newton-based ESC algorithm \cite{LIUFinal} for mitigating the effects of time delays.}
\label{Fig1}
\end{figure}

A constant and known delay matrix $D \ge 0$ is introduced in the actuation path, so that the measured output is given by (see Figure  \ref{Fig1}):
\begin{equation} \label{eq:1}
    y(t) = Q(\textcolor{black}{\theta^D(t)}). 
\end{equation}

The results found in this paper can be directly extended to the measurement delay case, as any equal delays in the input channels can be redirected to the output of the static map. The case where input delays 
and output delays 
occur simultaneously can also be addressed by assuming that the total delay to be neutralized is simply the sum of the parts. 
In a more general framework, we consider the following representation of the input-output delay:
  \begin{equation} \label{eq:AtrasoEntr}
   \textcolor{black}{\theta^D(t):=\theta(t-D)} = e^{-Ds}[\theta(t)],
\end{equation}
where it is assumed that the inputs have distinct, known, and ordered delays: 
\begin{equation} \label{eq:AtrasosDistintos}
    D=\textrm{diag}\{D_1, D_2, \cdot \cdot \cdot, D_n \}, \quad 0 \leq D_1 \leq D_2 \leq \ldots \leq D_n.
\end{equation}

In addition, we assume that the constants $D_i$ are known, $\forall i \in \{1,2,...,n \}$.

\subsection{System and Signals}
Without loss of generality, we assume a maximum seeking problem, where the goal is to find the maximizing value of $\theta \in \mathbb{R}^n$  denoted by $\theta^*\in \mathbb{R}^n$. For simplicity, we consider the quadratic nonlinear map:
\begin{equation} \label{eq:2}
y(\theta)=y^\ast+\frac{1}{2}(\theta-\theta^\ast)^TH(\theta-\theta^\ast),
\end{equation}
\noindent 
where the extremum point is $y^\ast \in \mathrm{R}$, and $H=H^T<0$ is the $n\times n$ unknown Hessian matrix of the corresponding static map. \textcolor{black}{Our methodology guarantees stability as long as the static map is at least locally quadratic, ensuring convergence within a neighborhood of the extremum. For maps that do not exhibit local quadratic behavior and may not yield exponential stability of the averaged system, the approach presented in \cite{ref_35_TAC,ref_389_HP_Krstic} could achieve asymptotic practical stability instead. However, their averaging framework does not apply to systems defined on Banach spaces, such as those governed by PDEs or delays.}

\textcolor{black}{
Without loss of generality, we assume that the input channels have distinct delays, ordered to satisfy (\ref{eq:AtrasosDistintos}). Any other ordering could be adopted, as the inputs can always be rearranged to comply with (\ref{eq:AtrasosDistintos}), which may result in a different Hessian matrix $H$ in (\ref{eq:2}).
}

Substituting (\ref{eq:2}) into (\ref{eq:1}), one obtains:
\begin{equation} \label{eq:Saida2}
    y(\theta)=y^\ast+\frac{1}{2}(\theta(t-D)-\theta^\ast)^TH(\theta(t-D)-\theta^\ast).
\end{equation}

From equation (\ref{eq:AtrasoEntr}), in the general case of multiple and distinct delays in the input channels, the delayed input vector can be represented as:
\begin{equation} \label{eq:EntradaTheta}
    \theta(t-D):= \begin{bmatrix}
        \theta_1(t-D_1) \\
        \theta_2(t-D_2) \\
        \vdots \\
        \theta_n(t-D_n)
    \end{bmatrix}.
\end{equation}

Let $\hat{\theta}$ be the estimate of $\theta^{\ast}$ and define the estimation error as
\begin{equation}\label{eq:3}
    \Tilde{\theta}(t):=\Hat{\theta}(t)-\theta^\ast.
\end{equation}
\noindent
\textcolor{black}{From the block diagram in Figure \ref{Fig1} and equation (\ref{eq:3}), we obtain $\dot{\tilde{\theta}}(t)=\dot{\hat{\theta}}(t)=U(t)$, since $\theta^*$ remains constant. By applying a delay of $D_i$ to both sides of this equation for the corresponding state variables, we derive the following expression for the error dynamics:}
\begin{equation} \label{eq:4}
    \dot{\tilde{\theta}}(t-D)=\begin{bmatrix}
U_1(t-D_1) \\
U_2(t-D_2) \\
\vdots      \\
U_n(t-D_n)
\end{bmatrix}\,, \;\;\; \dot{\tilde{\theta}}_i(t-D_i)=U_i(t-D_i),
\end{equation}
\noindent
with the average version given by:
\begin{equation} \label{eq:5b}
    \dot{\tilde{\theta}}_i^{av}(t-D_i)=U_i^{av}(t-D_i).
\end{equation}

From Figure \ref{Fig1}, we also have: 
\begin{equation} \label{eq:5}
    G(t)=M(\eta(t-D_i))y(t)
    ~~~~\text{and}~~~~ \theta(t)=\Hat{\theta}(t)+S(\eta(t)).
\end{equation}

The stochastic sinusoidal perturbation signals, also referred to as dither signals, $S(\eta(t)) \in \mathbb{R}^n$ and $M(\eta(t)) \in \mathbb{R}^n$ are defined as follows::
\begin{equation} \label{eq:6}
    S(\eta(t))=\bigg[a_1\sin(\eta_1(t)) \: ... \: a_n\sin(\eta_n (t)) \bigg]^T, 
\end{equation}
\begin{equation} \label{eq:7}
    M(\eta(t))=\bigg[\frac{2}{a_1}\sin(\eta_1 (t)) \: ... \: \frac{2}{a_n}\sin(\eta_n (t))\bigg]^T,
\end{equation}
\noindent
where $a_i$ are nonzero amplitudes $\forall i=1,2,...,n$. The subscript $i$ denotes the $i$-th entry of the vector $\eta_i(t)$. The elements of the stochastic Gaussian perturbation vector $\eta(t)$ are sequential and mutually independent, such that  $\mathbb{E}\{\eta(t)\}=0, \mathbb{E}\{\eta^2_i(t)\}=\sigma^2_i$ and $\mathbb{E}\{\eta_i(t)\;\eta_j(t)\}=0$, $\forall i \neq j$, with $\mathbb{E}\{\cdot\}$ denoting the expected value of the signal. Moreover, we assume that the probability density function of the perturbation vector is symmetric about its average.

Stochastic sinusoidal perturbations are generated using the
standard Brownian motion process $W_{\omega t}$ (also known as the Wiener process) over the boundary of a circle \citep{c6,c14}. Accordingly, we refer to this Markov Process as a stochastic process that satisfies the Markov property with respect to the natural filtration. 

From Figure \ref{Fig1}, the input vector signal is constructed as follows: 
\begin{equation} \label{eq:8}
    \theta(t) = \hat{\theta}(t)+S(\eta(t)),
\end{equation}
where $\hat{\theta}$ represents the estimate of $\theta^\ast$ and 
\begin{equation} \label{eq:9}
    \eta_i(t)= \omega \pi(1 + \text{sin}(W^i_{\omega t}))
\end{equation}
\noindent
defines a homogenous ergodic Markov process with a nonzero frequency $\omega$. The terms $W^i_{\omega t}$ indicate that different Wiener processes are mutually independent for each channel and for $\omega >0$. Using the time scale $\tau=\omega t$ and applying the stochastic chain rule \citep{c6}, we obtain: 
\begin{equation} \label{eq:10}
    d\eta_i =-\frac{\pi}{2}\text{sin}(W^i_{\tau})d\tau + \pi \text{cos}(W^i_{\tau})dW^i_{\tau}.
\end{equation}

The estimate of the unknown Hessian $H$ is given by 
\begin{equation} \label{eq:11}
    \hat{H}(t)=N(\eta(t-D))y(t),
\end{equation}
which satisfy the following average property \footnote{The ergodic property implies that the time average of a function of the process along its trajectories exists almost surely and equals the space average: $\lim_{T \to \infty}\frac{1}{T}\int_0^{T}a(Z_t)dt =\int_{S_Z}a(z)\mu(dz)$,  for any integrable function $a(\cdot)$, where $\mu(dz)$ is the invariant distribution of the samples $Z_t$ over $S_Z$ \cite{c14}.}
\begin{equation} \label{eq:12}
    \frac{1}{\Pi}\int\limits_0^\Pi N(\sigma)yd\sigma=H, \quad \Pi=2\pi/\omega.
\end{equation}

Equation (\ref{eq:12}) was demonstrated in \citep{c4,LIUFinal} for the case of a quadratic nonlinear map as in (\ref{eq:2}). Therefore, the average version of $H(t)$ is given by $\hat{H}_{av}(t)=(Ny)_{av}(t)=H$.

The elements of the $n \times n$ demodulation matrix $N(\eta(t))$, used to obtain the Hessian estimate, are given by:
\begin{equation} \label{eq:13}
    N_{ij}(t) = \left \{ \begin{matrix} \frac{16}{a_i^2}\Bigg(\text{sin}^2(\eta_i(t-D_i))-\frac{1}{2}\Bigg), & i=j \\ \frac{4}{a_ia_j}\text{sin}(\eta_i(t-D_i))\text{sin}(\eta_j(t-D_j)), & i \neq j. \end{matrix} \right.
\end{equation}

Thus, following  \citep{c4,LIUFinal}, we define the measurable signal:
\begin{equation} \label{eq:14}
    z(t)=\Gamma(t)G(t),
\end{equation}
and, by employing  averaging computation, we verify from equations (\ref{eq:5}), (\ref{eq:6}) and (\ref{eq:14}) that:
\begin{equation} \label{eq:14b}
    z_{av}(t)=\frac{1}{\Pi}\int_0^\Pi \Gamma M(\lambda)yd \lambda = \Gamma_{av}(t)H\Tilde{\theta}_{av}(t-D),
\end{equation}
where $\Gamma(t)$ is governed by the following Riccati differential equation: 
\begin{equation} \label{eq:15}
    \dot{\Gamma}=\omega_r\Gamma-\omega_r\Gamma\Hat{H}\Gamma,
\end{equation}
\noindent
where $\omega_r>0$ is a design constant. Equation (\ref{eq:15}) provides an estimate of the Hessian's inverse ($H^{-1}$), avoiding direct inversion of the estimated Hessian, which may assume null values during the transient phase. The estimation error of the Hessian's inverse is defined as:
\begin{equation} \label{eq:16}
    \tilde{\Gamma}(t)=\Gamma(t)-H^{-1}
\end{equation}
\noindent
with its dynamics given by the following equation, derived from (\ref{eq:15}) and (\ref{eq:16}):
\begin{equation} \label{eq:17}
   \dot{\tilde{\Gamma}} = \omega_r \left [ \tilde{\Gamma}+H^{-1} \right ] \times \left [ I_{n\times n}-\Hat{H}(\tilde{\Gamma}+H^{-1}) \right ],
\end{equation}
where $I_{n\times n}$ denotes the identity matrix of order $n$.

\subsection{Predictor Feedback via Hessian's Inverse Estimation}
By using averaging analysis \citep{c1,c6}, we can verify from $G(t)$ in (\ref{eq:5}) and $z(t)$ in (\ref{eq:14}) that: 
\begin{equation} \label{eq:21}
    z_{av}(t)= \Gamma_{av}(t)H\tilde{\theta}_{av}(t-D).
\end{equation}

Using (\ref{eq:16}), equation (\ref{eq:21}) can be rewritten in terms of $\tilde{\Gamma}_{av}(t)=\Gamma_{av}(t)-H^{-1}$ as: 
\begin{equation} \label{eq:22}
    z_{av}(t)=\tilde{\theta}_{av}(t-D)+\tilde{\Gamma}_{av}(t)H\tilde{\theta}_{av}(t-D).
\end{equation}

The second term on the right-hand side of (\ref{eq:22}) is quadratic in $(\tilde{\Gamma}_{av},\tilde{\theta}_{av})$; thus, linearizing $\Gamma_{av}$(t) around $H^{-1}$ yields a linearized form: 
\begin{equation} \label{eq:23}
    z_{av}(t)=\tilde{\theta}_{av}(t-D),
\end{equation}
\noindent
where $\Tilde{\theta}_{av}(t-D)=[\Tilde{\theta}_1^{av}(t-D_1) ~~...~~ \Tilde{\theta}_n^{av}(t-D_n]^T$. Therefore, from (\ref{eq:5b}) and (\ref{eq:23}), we obtain the following form of the average model with state $z_{av}(t)$: 
\begin{equation} \label{eq:23b}
    \Dot{z}_{av}(t)=\begin{bmatrix}
        \Dot{z}_1^{av}(t) \\
        \Dot{z}_2^{av}(t) \\
        \vdots \\
        \Dot{z}_n^{av}(t)
    \end{bmatrix} = \begin{bmatrix}
        U_1^{av}(t-D_1) \\
        U_2^{av}(t-D_2) \\
        \vdots \\
        U_n^{av}(t-D_n)
    \end{bmatrix}.
\end{equation}

From (\ref{eq:4}) and  (\ref{eq:23}), we obtain the following average model:
\begin{equation} \label{eq:24}
    \dot{\tilde{\theta}}_{av}(t-D)=U_{av}(t-D),
\end{equation}
\begin{equation} \label{eq:25}
    \dot{z}_{av}(t)=U_{av}(t-D),
\end{equation}
\noindent
where  $U_{av} \in \mathbb{R}^{n}$ is the average control derived  from $U \in \mathbb{R}^{n}$.

To motivate the design of predictor feedback, the main idea here is to compensate for the delay by incorporating the future state $z(t+D)$ into the feedback control law, or $z_{av}(t+D)$ into the system’s average model. 
Using the variation of constants formula for equations (\ref{eq:24}) and (\ref{eq:25}), the future state is given by:
\begin{equation} \label{eq:26}
    z_{av}(t+D)=z_{av}(t)+\int\limits_{t-D}^t U_{av}(\sigma)d\sigma,
\end{equation}
where the integral of $U_{av}(\sigma)$ in (\ref{eq:26}) is computed  over the past interval $[t-D,t]$. Given any diagonal matrix $K>0$ as a stabilizing gain, the average control must be:
\begin{equation} \label{eq:27}
    U_{av}(t)= -K \left [z_{av}(t) +\int\limits_{t-D_i}^t U_{av}(\sigma)d\sigma  \right ],
\end{equation}
\noindent
which yields the average control $U_{av}(t)=-Kz_{av}(t+D), \forall{t}>0$, as desired. 
Therefore, the closed-loop average system can be expressed by:
\begin{equation} \label{eq:28}
 \frac{d\tilde{\theta}_{av}(t)}{dt}=-K\tilde{\theta}_{av}(t),  \quad \forall{t} \geq D_n,
\end{equation}
\textcolor{black}{since the quadratic term $\tilde{\Gamma}_{av}H\tilde{\theta}_{av}$ in (\ref{eq:22}) was 
neglected in the linearization of the system (\ref{eq:28}).} Thus, the Newton algorithm achieves local exponential stability with a guaranteed convergence rate determined by $-K$, being independent of the unknown Hessian $H$.

In \citep{c7,Oliveira2022}, it is shown that the control objectives can be achieved by applying a simple modification to the predictor-based controller, specifically incorporating a low-pass filter, so that the average theorem in infinite dimensions \citep{AHale,ALehman} can be invoked. In this case, the following averaging-based predictor feedback is proposed to compensate for the time delays:
\begin{equation} \label{eq:29}
    U_i(t)=\frac{c_i}{s+c_i}\left \{ -k_i  \left [ z_i(t)+\int\limits_{t-D_i}^t U_i(\tau)d\tau  \right ]         \right \},
\end{equation}
\noindent
for all $k_i$, $c_i>0$, $i=1,2,...,n$, and $c_i$ sufficiently large. In other words, the predictor feedback takes the form of a low-pass filter applied to a non average version of (\ref{eq:27}). Using a mixed notation of time and frequency domains in (\ref{eq:29}), the transfer function acts as an operator on a time-domain function.  The predictor feedback (\ref{eq:29}) is average based because $z(t)$ in (\ref{eq:14}) is updated according to the estimate $\Gamma(t)$ of the unknown Hessian inverse $H^{-1}$, as given by (\ref{eq:15}), with $\Hat{H}(t)$ in (\ref{eq:11}) satisfying the average property (\ref{eq:12}).

\section{Stability Analysis}\label{sec4}
The following theorem summarizes the stability and convergence properties of the closed-loop feedback system. The operators $\mathbb{E}\left \{ \cdot \right \}$ and $\mathbb{P}\left \{ \cdot \right \}$ denote the expected value and probability of the corresponding signals, respectively. 

\textbf{Theorem 1} Consider the closed-loop system in Figure \ref{Fig1}, with multiple distinct delays (\ref{eq:AtrasosDistintos}) applied to the inputs of the locally quadratic map (\ref{eq:2}). There exists $c^{\ast}>0$ such that, $\forall{c_i} \geq c^{\ast}$, $\exists$ $\omega^{\ast}(c_i)$ 
such that $\forall$ $\omega \geq \omega^{\ast}$, the delayed closed-loop system described by (\ref{eq:4}) and (\ref{eq:29}), with $z(t)$ in (\ref{eq:14}), $G(t)$ in (\ref{eq:5}), $\Gamma(t)$ in (\ref{eq:15}) and state $\tilde{\Gamma}(t)$, $\tilde{\theta_i}(t-D_i)$, and $U_i(\sigma)$ $\forall{\sigma} \in [t-D_i,t]$, and $\forall i \in \{1,2,...,n\}$, has a unique exponentially stable solution satisfying:
\begin{equation} \label{eq:30}
\mathbb{E} \Bigg\{  |\tilde{\Gamma}(t)|^2 + \sum_{i=1}^n\Big[\tilde{\theta}_i(t-D_i)\Big]^2  + \left [ U_i(t) \right]^2  
     +\int\limits_{t-D}^t \left [ U_i(\sigma) \right]^2 d\sigma \Bigg\}^{1/2}\le \mathcal{O}(1/\omega), \quad \forall{t} \rightarrow +\infty.
\end{equation} 
In particular,
\begin{equation} \label{eq:31}
  \lim_{(1/\omega) \to 0}\mathbb{P} \left \{ \lim _{t \to \infty} sup\:|\theta(t) - \theta^{\ast}| \le O(|a|+1/\omega) \right \} = 1,
\end{equation}
\begin{equation} \label{eq:32}
  \lim_{(1/\omega) \to 0}\mathbb{P} \left \{ \lim _{t \to \infty} sup\:|y(t) - y^{\ast}| \le O(|a|^2+1/\omega^2) \right \} = 1.
\end{equation}
where $a=\big [a_1,\;a_2,\; ...\;,\; a_n \big]^T$.

The steps of the proof are presented in Sections $3.1$ to $3.8$ below.

\subsection{Transport PDE for Delay Representation}
According to \citep{c1}, the delay in (\ref{eq:4}) can be represented using  a transport PDE as follows:
\begin{equation} \label{eq:33}
    \dot{\tilde{\theta_i}}(t-D_i)=u_i(0,t),
\end{equation}
\begin{equation} \label{eq:34}
    \partial_t u_i(x,t)=\partial_x u_i(x,t), \quad x \in [0,D_i],
\end{equation}
\begin{equation} \label{eq:35}
    u_i(D_i,t)=U_i(t),
\end{equation}
where the solution to (\ref{eq:34})-(\ref{eq:35}) is:
\begin{equation} \label{eq:36}
    u_i(x,t)=U_i(t+x-D_i)\;,
\end{equation} 
where $t$ is the time variable, $D_i$ is the delay, and $x$ is the spatial variable of the PDE.

\subsection{Equations of the Closed-Loop System}
Substituting (\ref{eq:5}) and (\ref{eq:11}) into (\ref{eq:29}) and representing the integrand in (\ref{eq:29}) by the transport PDE state, we obtain:
\begin{equation} \label{eq:38}
    U_i(t)=\frac{c_i}{s+c_i}\left \{-k_i \Bigg[z_i(t) + \int\limits_{0}^{D_i} u_i(\tau,t)d\tau \Bigg] \right\}.
\end{equation}

Finally, by substituting (\ref{eq:38}) into (\ref{eq:35}), we can rewrite (\ref{eq:33})-(\ref{eq:35}) as: 
\begin{equation} \label{eq:39}
    \dot{\tilde{\theta_i}}(t-D_i)=u_i(0,t),
\end{equation}
\begin{equation} \label{eq:40}
    \partial_t u_i(x,t)=\partial_x \;u_i(x,t)\;, \quad x \in [0,D_i],
\end{equation}
\begin{equation} \label{eq:41}
    u_i(D_i,t)=\frac{c_i}{s+c_i}\left \{-k_i z(t) - k_i\int\limits_{0}^D u (\tau,t)d\tau \right\}.
\end{equation}

\subsection{Average Model of the Closed-loop System}
Now, defining
\begin{equation} \label{eq:43}
    \tilde{\vartheta}(t):=\tilde{\theta}(t-D)
\end{equation}
and using (\ref{eq:23}) and (\ref{eq:29}), the average version of the system (\ref{eq:39})-(\ref{eq:41}) becomes: 
\begin{equation} \label{eq:44}
    \dot{\tilde{\vartheta}}_i^{av}(t)=u_i^{av}(0,t),
\end{equation}
\begin{equation} \label{eq:45}
    \partial_t u_i^{av}(x,t)=\partial_x \; u_i^{av}(x,t)\;, \quad x \in [0,D_i]\;,
\end{equation}
\begin{equation} \label{eq:46}
    \frac{d}{dt}u_i^{av}(D_i,t)=-c_i u_i^{av}(D,t) - c_ik_i \Bigg [\tilde{\vartheta}_i^{av}(t)+\int\limits_0^Du_i^{av}(\sigma,t)d\sigma \Bigg],  
\end{equation}
\noindent
where the filter $c_i/(s + c_i)$ is represented in state-space form. Meanwhile, the average model for the estimation error associated with the Hessian's inverse in (\ref{eq:17}) is $\frac{d\tilde{\Gamma}_{av}(t)}{dt}=-\omega_r\tilde{\Gamma}_{av}(t)-\omega_r\tilde{\Gamma}_{av}(t)H\tilde{\Gamma}_{av}(t)$ with its linearized form given by: 
\begin{equation} \label{eq:48}
    \frac{d\tilde{\Gamma}_{av}(t)}{dt}= - \omega_r\tilde{\Gamma}_{av}(t).
\end{equation}

\subsection{Backstepping Transformation, Its Inverse, and the Target System}
Consider the following infinite-dimensional backstepping transformation \citep{c1} of the delayed state:
\begin{equation} \label{eq:49}
    w_i(x,t)=u_i^{av}(x,t) + k_i \left [ \tilde{\vartheta}_i^{av}(t) + \int\limits_{0}^x u_i^{av} (\sigma,t)d\sigma\right],
\end{equation}
\noindent
where $\vartheta_{av}(t):=z_{av}(t)=\tilde{\theta}_{av}(t-D)$, as given by (\ref{eq:23}). The transformation in   (\ref{eq:49}) maps the linearized system (\ref{eq:44})-(\ref{eq:46}) to:
\begin{equation} \label{eq:50} 
    \dot{\tilde{\vartheta}}_i^{av}(t)=-k_i \tilde{\vartheta}_i^{av}(t) +w_i(0,t),
\end{equation}
\begin{equation} \label{eq:51}
    \partial_t w_i(x,t)=\partial_x w_i(x,t), \quad x \in [0,D_i],
\end{equation}
\begin{equation} \label{eq:52}
    w_i(D_i,t)=-\frac{1}{c_i}\partial_t u_i^{av}(D_i,t), \quad i=1,2,...,n.
\end{equation}
We observe that (\ref{eq:50}) is input-to-state stable \cite{Karafyllis2018} with respect to $w_i(0,t)$ and (\ref{eq:51})-(\ref{eq:52}) exhibit finite-time stability as $c_i\to +\infty$, implying $w_i(D_i,t) \to 0$.

By using  (\ref{eq:39}), we partially differentiate the transformed state $w_i(x,t)$ in (\ref{eq:49}) with respect to time $t$ and consider $x=D_i$, yielding:
\begin{equation} \label{eq:53}
    \partial_t w_i(D_i,t)=\partial_tu_i^{av}(D_i,t)+k_iu_i^{av}(D_i,t),
\end{equation}
where $\partial_tu_i^{av}(D_i,t)=\dot{U}_i^{av}(t)$. Furthermore, the transformation in (\ref{eq:49}) is invertible, given by:
\begin{equation} \label{eq:54}
    u_i^{av}(x,t)=w_i(x,t)-k_i \Bigg [e^{-k_ix}\tilde{\vartheta}_i^{av}(t)
    +\int\limits_{0}^x e^{-k_i(x-\sigma)} w_i(\sigma,t)d\sigma \Bigg ].
\end{equation}

\subsection{Lyapunov-Krasovskii Functional}
The exponential stability of the overall system is achieved with the following Lyapunov-Krasovskii functional:
\begin{equation} \label{eq:Lya0}
    V(t)= \sum_{i=1}^n V_i(t), \quad i=1,...,n,
\end{equation}
\noindent
where $V_i(t)$ are functionals defined as:
\begin{align} \label{eq:Lya1}
V_i(t) = \frac{1}{2}\Big[\tilde{\vartheta}_i^{av}(t)\Big]+\frac{\overline{a}_i}{2} \int_0^{D_i}(1+x)w_i(x,t)^2dx   +\frac{1}{2}w_i^2(D_i,t), 
\end{align}
for each subsystem in (\ref{eq:50})-(\ref{eq:52}), with   $\overline{a}_i$ to be chosen later. Thus, we have: 
\begin{align} \label{eq:Lya2}
&\dot{V}_i(t)= -k_i\Big[\tilde{\vartheta}_i^{av}(t)\Big]^2+\tilde{\vartheta}_i^{av}(t)\omega_i(0,t)
+ \overline{a}_i\int\limits_{0}^{D_i}(1+x)\omega_i(x,t) \partial_x \omega_i(x,t) \nonumber \\
& +\omega_i(D_i,t) \partial_t\omega_i(D_i,t).
\end{align}
\noindent
By employing Young's inequality and integrating  by parts, we obtain: 
\begin{align} \label{eq:Lya3}
&    \dot{V}_i(t) \le -k_i \Big[\tilde{\vartheta}_i^{av}(t)\Big]^2 + \frac{\Big[\tilde{\vartheta}_i^{av}(t)\Big]^2}{2\overline{a}_i}-\frac{\overline{a}_i}{2}\int\limits_{0}^{D_i}\omega_i^2(x,t)dx \nonumber \\
&  +\omega_i(D_i,t) \Bigg[\partial_t \omega_i(D_i,t)+\frac{\overline{a}_i(1+D_i)}{2}\omega(D_i,t) \Bigg].
\end{align}
\noindent 
Recalling that $k_i>0$, we choose $\overline{a}_i=1/k_i$. Hence, 
\begin{equation} \label{eq:Lya4}
    \dot{V}_i(t)=-\frac{1}{2\overline{a}_i}\Big[\tilde{\vartheta}_i^{av}(t)\Big]^2-\frac{\overline{a}_i}{2}\int\limits_{0}^{D}\omega_i^2(x,t)dx +\omega_i(D_i,t)\bigg [\partial_t \omega_i(D_i,t)+\frac{\overline{a}_i(1+D_i)}{2}\omega_i(D_i,t) \Bigg ].
\end{equation}

Now, considering (\ref{eq:Lya4}) along with (\ref{eq:54}) and completing the squares, we have:
\begin{align} \label{eq:Lya5}
&    \Dot{V}_i(t)=-\frac{1}{4\overline{a}_i}-\frac{\overline{a}_i}{4}\int^{D_i}\omega^2_i(x,t)dx+\overline{a}_i\Big|k_i^2 e^{-k_iD_i}\Big|^2\omega_i^2(D_i,t) \nonumber \\
& +\frac{1}{\overline{a}_i}\Big|\Big|k_i^2e^{-k_i(D_i-\sigma)}\Big|\Big|\omega_i^2(D_i,t) +\Bigg[\frac{\overline{a}_i(1+D_i)}{2}+k_i\Bigg]\omega_i^2(D_i,t) - c_i\omega_i^2(D_i,t),
\end{align}
where we have applied the Cauchy-Schwartz and Young
inequalities to obtain  the following expressions: 
\begin{align} \label{eq:64}
& -w_i(D_i,t)^T\Big\langle(KH)^2e^{KH(D_i-\sigma)},w_i(\sigma,t)\Big\rangle \nonumber \\ &\le \big|w_i(D_i,t)\big| \Big| \Big|(KH)^2e^{KH(D_i-\sigma)}\Big|\Big|\big|\big|w_i(t)\big|\big| \nonumber \\
& \le \frac{a}{4}\big|\big|w_i(t)\big|\big|^2+\frac{1}{a}\Big|\Big|(KH)^2e^{KH(D_i-\sigma)}\Big|\Big|^2w_i^2(D_i,t).
\end{align}
\noindent

The notation $\langle \cdot, \cdot \rangle$ denotes the inner product over the spatial variable $\sigma \in [0,D_i]$, upon which both $e^{KH(D_i-\sigma)}$ and $w_i(\sigma,t)$ depend. Additionally, $||\cdot||$ represents the $\mathcal{L}_{2}[0,D_i]$ norm in $\sigma$. Therefore, from (\ref{eq:Lya5}), we can derive the following inequality:
\begin{equation} \label{eq:Fim0}
    \Dot{V}_i(t) \le -\frac{1}{4\overline{a}_i}\Big[ \Tilde{\vartheta}_i^{av}(t)\Big] - \frac{\overline{a}_i}{4}\int_0^{D_i}\omega^2_i(x,t)dx-(c_i-c_i^\ast)\omega_i^2(D_i,t),
\end{equation}
\noindent
where
\begin{equation} \label{eq:Fim1}
    c_i^\ast=\frac{\overline{a}_i(1+D_i)}{2}+k_i+\overline{a}_i\Big|k_i^2e^{k_iD_i}\Big|+\frac{1}{\overline{a}_i}\Bigg|\Bigg|k_i^2e^{-k_i(D_i-\sigma)}\Bigg|\Bigg|.
\end{equation}

Hence, from (\ref{eq:Fim0}), if $c_i$ is chosen such that $c_i>c_i^\ast$, we obtain:
\begin{equation}\label{eq:LyaFinal}
    \dot{V}_i(t)\le - \mu_i V_i(t) \quad \textrm{or} \quad \dot{V}(t)\le - \mu V(t),
\end{equation}
for some $\mu_i>0$ and $\mu=\min(\mu_i)$. Therefore, the closed-loop system is exponentially stable in the sense of the full-state norm
\begin{equation} \label{eq:LyaFinal2}
    \sqrt{|\tilde{\vartheta}_{av}(t)|^2+\sum_{i=1}^nw_i^2(D_i,t)+\sum_{i=1}^n\int\limits_{0}^{D_i}w_i^2(x,t)dx},
\end{equation}
i.e., in the transformed variable $(\tilde{\vartheta}_{av},w)$.

\subsection{Exponential Stability for the Average System in the Original State Variables}
The objective here is to establish the exponential stability of the closed-loop average systems in the terms of the norm
\begin{equation} \label{eq:56}
    \Upsilon(t)=\sqrt{|\tilde{\vartheta}_{av}(t)|^{2}+\sum_{i=1}^n[u_i^{av}(D_i,t)]^2 +\sum_{i=1}^n\int\limits_{0}^{D_i} (u_i^{av}(x,t))^2dx}.
\end{equation}
This result is achieved by using equations (\ref{eq:49}), (\ref{eq:54}),  (\ref{eq:Lya0}), (\ref{eq:Lya1}), and noting the existence of  positive constantes $\alpha_1$ and $\alpha_2$ such that $\alpha_1\Upsilon(t)^2 \le V(t) \le \alpha_2 \Upsilon(t)^2$. Applying the Cauchy-Schwartz inequality along with additional calculations, as outlined in \citet[Chapter~2]{c1} and using  (\ref{eq:LyaFinal}), we obtain:
\begin{equation} \label{eq:153}
    \Upsilon(t) \le \frac{\alpha_2}{\alpha_1}e^{-\mu t}\Upsilon(0),
\end{equation}
\noindent
which completes the proof of exponential stability in the original variable $(\Tilde{\vartheta}_{av},u_{av})$.

\subsection{Invoking the Stochastic Averaging Theorem}
From (\ref{eq:4}), (\ref{eq:17}) and (\ref{eq:29}), we have:
\begin{equation} \label{eq:60}
\begin{split}
\hspace{-0.5cm}\frac{d}{dt}
\begin{bmatrix}
\tilde{\theta}_i(t-D_i)     \\
U_i(t) \\
\tilde{\Gamma}(t)   
\end{bmatrix} = 
\begin{bmatrix}
0 \\
-c_iU_i(t) \\
0   
\end{bmatrix} 
+\begin{bmatrix}
U_i(t-D_i)     \\
-c_i k_i z_i(t) -c_ik_i \int\limits_{t-D_i}^t U_i(\tau)d(\tau) \\
\omega_r\Big[\tilde{\Gamma}(t)+H^{-1}\Big]\times\Big[I_{n\times n}-\hat{H}(t)\Big(\tilde{\Gamma}(t)+H^{-1}\Big)\Big]
\end{bmatrix}.
\end{split}
\end{equation}

Now, we define the state vector
\begin{equation} \label{eq:61}
    \textbf{u}^\epsilon(t):=\begin{bmatrix}
\tilde{\theta}(t-D)     \\
U(t) \\
\tilde{\Gamma}(t)   
\end{bmatrix}.
\end{equation}
This definition allows us to express (\ref{eq:60}) in the form of a stochastic functional differential equation:
\begin{equation} \label{eq:62}
    \frac{d}{dt}\textbf{u}^{\epsilon}(t)=G(\textbf{u}_t^\epsilon)+\epsilon F(t,\textbf{u}_t^\epsilon, \eta(t), \epsilon),
\end{equation}
where $\epsilon:=1/\omega$. Since $\eta(\tau)$ is a homogeneous ergodic Markov process with an invariant measure $\mu(d\eta)$ and undergoes the property of exponential ergodicity, we have $\textbf{u}_t^\epsilon(\delta)=\textbf{u}^\epsilon(t+\delta)$ for $-D_n\le \delta \le 0$. 

Furthermore, the function $G:\textbf{C}_3([-D_n,0]) \rightarrow{\mathbb{R}^3}$, along with the Lipschitz function $F:\mathbb{R}_+ \times \textbf{C}_3([-D_n,0])\times Y \times [0, 1] \to \mathbb{R}^3$, with $F(\tau,\;0,\;\eta,\; \epsilon) = 0$, are continuous mappings, where $\textbf{C}_3([-D_n,0])$ denotes the class of continuous vector functions over the interval $[-D_n,\;0]$. 

Thus, we can apply the Averaging Theorem by \citet{Kata} to conclude the exponential $p$-stability result (with $p=2$) of the original random system as $\epsilon \rightarrow 0$, thereby obtaining inequality (\ref{eq:30}). 

\subsection{Practical Convergence to the Extemum Point}
We define the stopping time \citep{c14b} as:
\begin{equation} \label{eq:63}
    \tau_{\epsilon}^{\Delta(\epsilon)}\!:=\!inf\left \{ \forall{t}\!\geq\! 0 :| \textbf{u}^{\epsilon}(t)| \!>\! {M}|\textbf{u}^{\epsilon}(0)|e^{-\lambda t} \!\!+\!\mathcal{O}(\epsilon)\right\},
\end{equation}
\noindent
representing the first time when the error vector norm no longer satisfies the exponential decay property. Let 
$M>0$ and $\lambda>0$ be constants, and let $T(\epsilon):(0,1) \rightarrow \mathbb{N}$ be a continuous function. Then, analogously to \citep{c14b}, the norm of the error vector $|\textbf{u}^{\epsilon}(t)|$ converges both \emph{almost surely}
(a.s.) and \emph{in probability} to a value below a residue error $\Delta(\epsilon)=\mathcal{O}(\epsilon)$:
\begin{equation} \label{eq:sm1}
    \lim_{\epsilon \to 0} \textrm{inf} \{ \forall{t}\geq  0 : |\textbf{u}^\epsilon(t)| > M |\textbf{u}^\epsilon(0)|e^{-\lambda t} + \Delta \} = \infty, \; \textrm{a.s.},
\end{equation}
\begin{equation} \label{eq:sm2}
 \lim_{\epsilon \to 0} \mathbb{P}\{|\textbf{u}^\epsilon(t)| \leq M|\textbf{u}^\epsilon(0)|e^{-\lambda t} + \Delta\,, \forall{t} \in [0,T(\epsilon)]\}=1,
\end{equation}
\noindent
with $\lim_{\epsilon \to 0}\textrm{T}(\epsilon) = \infty$. Equations (\ref{eq:sm1}) shows that $\tau_{\epsilon}^{\Delta (\epsilon)}$ tends to infinity as $\epsilon$ approaches zero. Similarly, in (\ref{eq:sm2}), the deterministic function $T(\epsilon)$ tends to infinity as $\epsilon \rightarrow 0$. Therefore, exponential convergence is sustained over an arbitrarily long time interval, and each component of the error vector converges to within $\Delta(\epsilon)=\mathcal{O}(\epsilon)$, specifically $\Tilde{\theta}(t)$. We can express this as $\lim_{\epsilon \to 0}\mathbb{P} \Big \{ \textrm{lim sup}_{t \to \infty}$ $|\tilde{\theta}(t)|\leq \mathcal{O}(\epsilon) \Big \}= 1$. 

From equations (\ref{eq:3}) and (\ref{eq:8}), we have:
\begin{equation} \label{eq:sm11}
\theta(t)- \theta^{\ast} =  \tilde{\theta}(t)+S(\eta(t)). 
\end{equation}

Considering that the first term on the right-hand side of (\ref{eq:sm11}) is ultimately of order  $\mathcal{O}(\epsilon)$, and the second term is of order $\mathcal{O}(|a|)$, we obtain equation (\ref{eq:31}). Finally, from (\ref{eq:2}) and (\ref{eq:31}), we derive (\ref{eq:32}). 

\section{Simulations}

\textcolor{black}{
The numerical simulation in this section is inspired by the \textit{source seeking problem} discussed in \cite[Section~IX]{c7}. In this scenario, multi-variable extremum seeking is used to locate the source of a signal --- such as a chemical, acoustic, or electromagnetic field --- whose concentration profile is described in (\ref{eq:Saida2}). The strength of the field decays with distance and attains a local maximum at $y^*$, with an unknown maximizer $\theta^*$. This optimization is performed without direct measurement of the position vector $\theta = (\theta_1, \theta_2)$,  relying solely on the measured scalar output $y$. Additionally, the vehicle’s two actuation paths, $\theta_1$ and $\theta_2$, are subject to distinct input delays, $D_1$ and $D_2$, respectively.}

To evaluate the proposed stochastic multivariable extremum seeking control with delay compensation, we examine a static quadratic map defined as: 
\begin{equation} \label{eq:100}
    Q(\theta)=5+\frac{1}{2}\Big(2(\theta_1)^2+4(\theta_2-1)^2+4\theta_1(\theta_2-1) \Big).
\end{equation}

\textcolor{black}{The delays are set to $D_{1}=50$ and $D_{2}=100$ for simulations involving distinct delays, while for simulations with equal delays, both are set to $D_{1}=D_{2}=100$}. The extremum point  is located at $\theta^\ast=[0,1]^T$ with $y^\ast=5$, and the unknown Hessian of the map is given by: 
\begin{equation} \label{eq:Hessiana}
    H=-\begin{pmatrix} 
        2 &\; 2 \\ 
        2 &\; 4 \\
        \end{pmatrix}.
\end{equation}

The numerical simulations of the predictor in (\ref{eq:29})  were performed wiht the following parameters: $a_1=a_2=0.22$, $c=20$, $\omega_r=0.007$, $\Gamma(0)=\begin{bmatrix}
        -1/100 &\; 0 \\ 
        0&\; -1/200 \\
        \end{bmatrix}$ and $K=0.005I_{2\times 2}$, where $I_{2\times 2}$ is the identity matrix of order $2$.

\subsection{\textcolor{black}{Performance Evaluation of Newton-Based ESC with Equal Input Delays}}\label{Sec4.2}
\textcolor{black}{In this section, we assess the performance of the Newton-based ESC when both actuation paths experience the same delay, specifically $D_1 = D_2 = 100$. Unlike the case in Section \ref{Sec4.1}, which consideres distinct delays for each input, this uniform delay configuration, as indicated in Figure \ref{Fig1} and formulated in equation (\ref{eq:100}), enables a systematic evaluation of the controller’s ability to compensate for time delays while ensuring convergence to the optimal point. Additionally, this scenario provides a basis for comparison with the case involving distinct delays (Section \ref{Sec4.1}).}
%
%
\begin{figure}[htb!]
    \centering
    \includegraphics[scale=0.54]{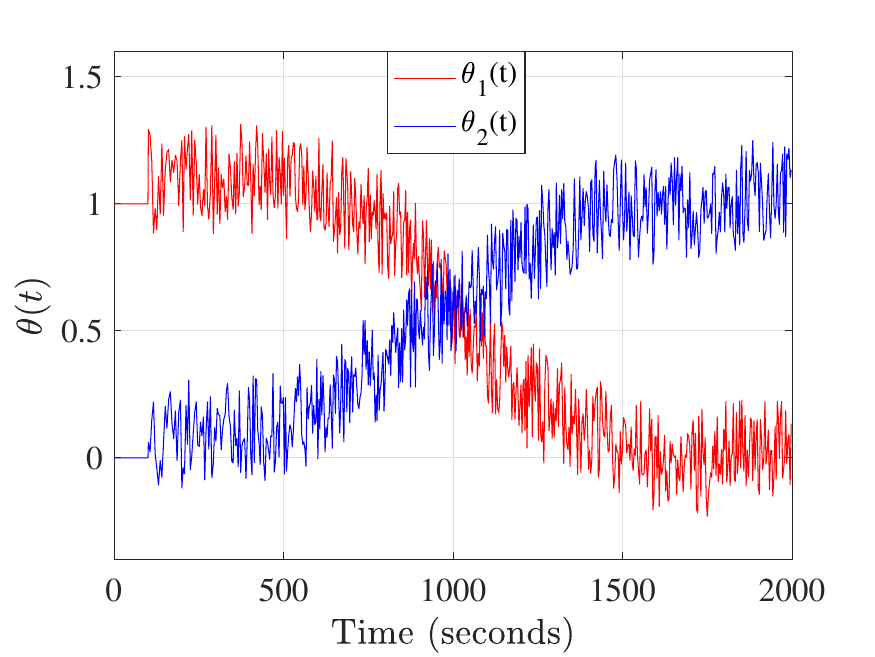}
    \caption{\textcolor{black}{Evolution of the system input $\theta(t)$ under the multivariable Newton-based stochastic ESC when both input channels experience equal delays.}}
    \label{fig:11}
\end{figure}


\textcolor{black}{Figure \ref{fig:11} shows the trajectory of the input $\theta(t)$ when the Newton-based ESC is applied with equal delays. Despite the presence of delay, the system remains stable, and the input converges to a region near the optimal point $\theta^* = (0,1)$. This indicates that the delay compensation mechanism integrated into the controller remains effective under this configuration.}



\subsection{\textcolor{black}{Performance Comparison: Newton-Based vs. Gradient-Based ESC}}\label{Sec4.3}

\textcolor{black}{In this subsection, we compare the performance of the gradient-based \cite{CODIT2024} and Newton-based stochastic ESC algorithms under different delay conditions. Both methods are tested according to the setup in equation (\ref{eq:100}), first in a delay-free environment and then in a scenario where both actuation paths experience equal delays, specifically $D_1 = D_2 = 100$. This comparison aims to assess their robustness and efficiency in handling time delays while ensuring convergence to the optimal point $y^*$.}

\begin{figure}[htb!]
    \centering
    \includegraphics[scale=0.54]{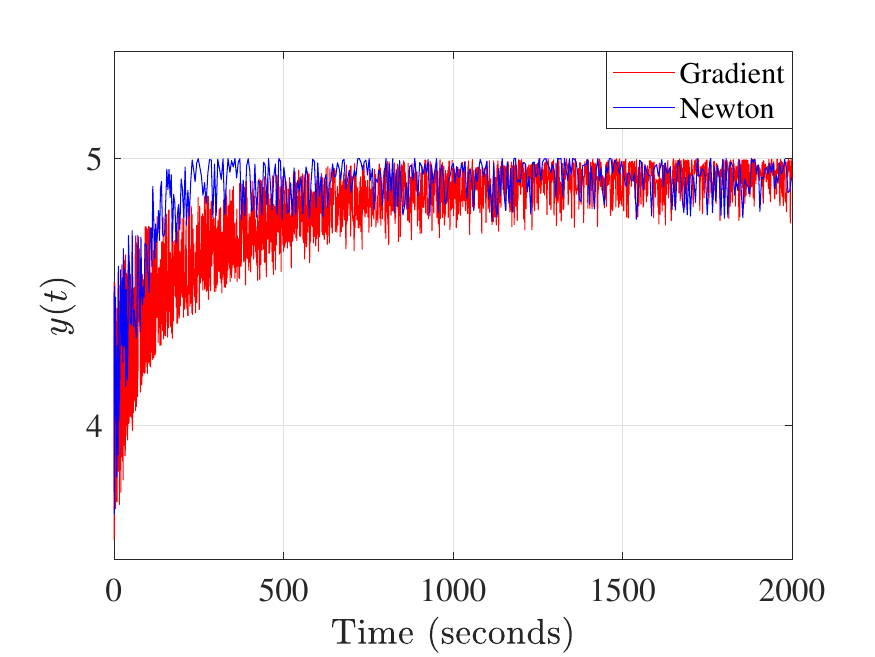}
    \caption{\textcolor{black}{Evolution of the system outputs $y(t)$ under the multivariable gradient-based and Newton-based stochastic ESC without delays.}}
    \label{fig:13}
\end{figure}

\begin{figure}[htb!]
    \centering
    \includegraphics[scale=0.54]{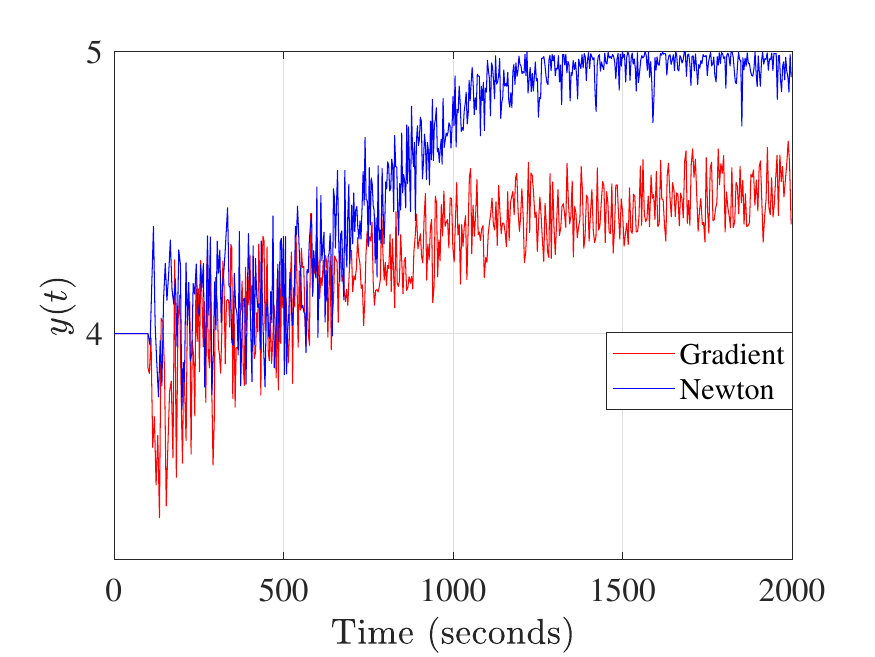}
    \caption{\textcolor{black}{Evolution of the system outputs $y(t)$ for the multivariable gradient-based and Newton-based stochastic ESC with predictor feedback when input channels experience equal delays.}}
    \label{fig:14}
\end{figure}

\textcolor{black}{Figure \ref{fig:13} presents the system output $y(t)$ for both ESC approaches in the absence of delays. The results show that while both controllers achieve convergence to the optimal point, the Newton-based ESC exhibits a significantly faster convergence rate compared to the gradient-based ESC. This highlights the advantage of using second-order information to accelerate optimization. Figure \ref{fig:14} extends this analysis to the delayed scenario, where equal delays are applied to both actuation paths. The results confirm that the Newton-based ESC maintains superior convergence properties even in the presence of delays, effectively compensating for time shifts and ensuring system stability. This evaluation underscores the advantages of Newton-based ESC, particularly in delay-sensitive applications.}



\subsection{\textcolor{black}{Performance Evaluation of Newton-Based ESC with Distinct Input Delays}}\label{Sec4.1}

\textcolor{black}{This section evaluates the performance of the proposed Newton-based stochastic ESC under distinct input delays. Unlike conventional methods that assume uniform delays, our approach compensates for each input delay, ensuring stability and convergence to the optimal point despite such distinct delay-induced disturbances.}

\begin{figure}[htb!]
    \centering
    \includegraphics[scale=0.54]{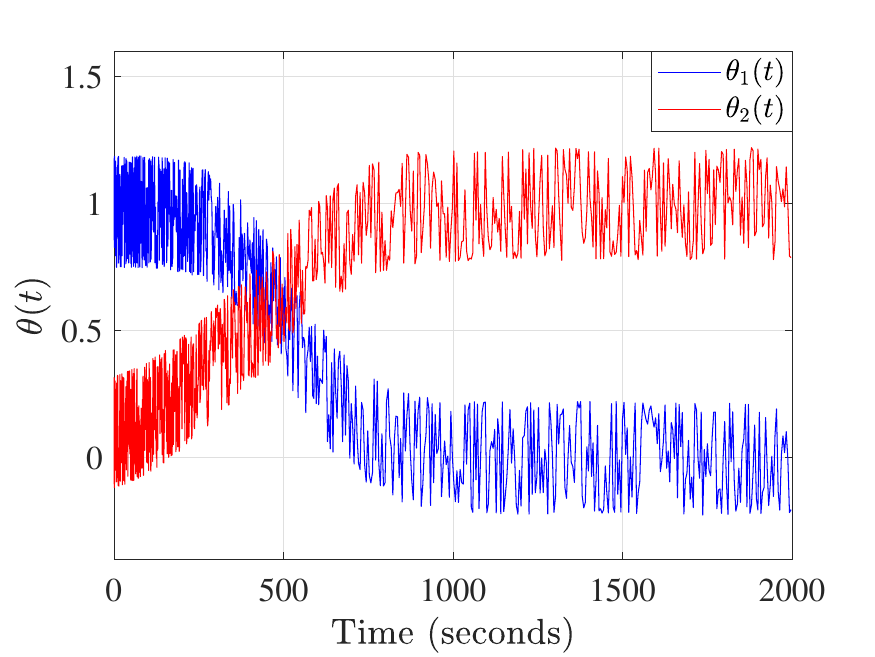}
    \caption{Evolution of the system input $\theta(t)$ under the multivariable Newton-based stochastic ESC without delays.}
    \label{fig:2}
\end{figure}

\begin{figure}[htb!]
    \centering
    \includegraphics[scale=0.54]{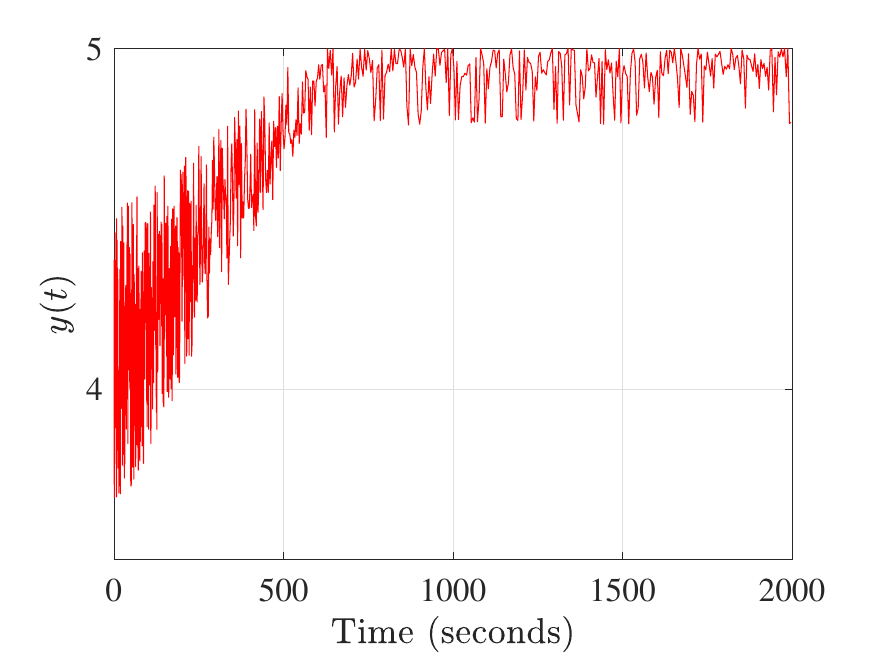}
    \caption{Evolution of the system output $y(t)$ under the multivariable Newton-based stochastic ESC without delays.}
    \label{fig:3}
\end{figure}

Figure \ref{fig:2} shows 
the convergence of the system input to a close vicinity around the optimal point  $\theta^\ast=[0,1]^T$.
Additionally, Figure \ref{fig:3} shows that the output of the nonlinear map converges to a small neighborhood around the optimal value $y^\ast=5$. In both figures, no delays were considered yet.

\begin{figure}[htb!]
    \centering    \includegraphics[scale=0.54]{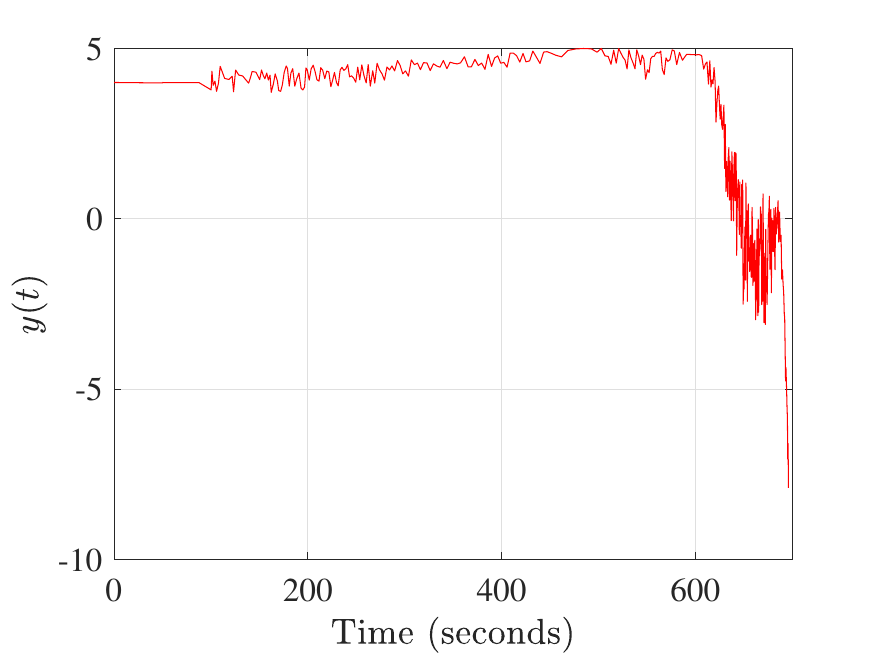}
    \caption{Evolution of the system output $y(t)$ under the multivariable Newton-based stochastic ESC without predictor feedback when input channels experience distinct delays.}
    \label{fig:4}
\end{figure}

In Figure \ref{fig:4}, we observe a phenomenon noted in the literature \citep{c7}: ESC is not robust in the presence of delays, causing the closed-loop system to become unstable when the delays are ignored and not properly compensated.

\begin{figure}[htb!]
    \centering    \includegraphics[scale=0.54]{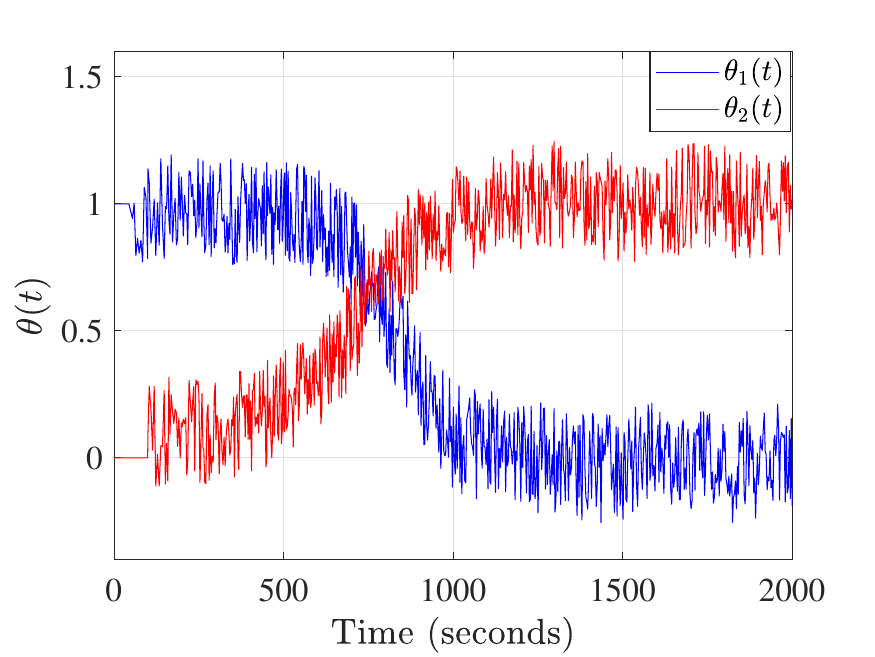} \caption{Evolution of the system input $\theta(t)$ under the multivariable Newton-based stochastic ESC and predictor feedback when input channels experience distinct delays.}
    \label{fig:5}
\end{figure} 

\begin{figure}[htb!]
    \centering
    \includegraphics[scale=0.54]{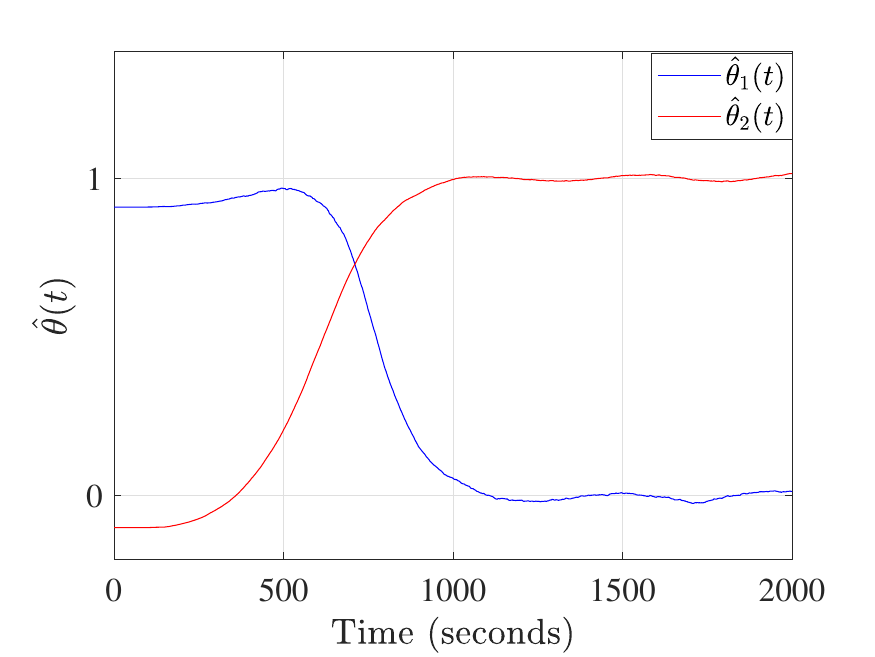}
    \caption{Estimate $\hat{\theta}(t)$ for the multivariable Newton-based stochastic ESC and predictor feedback with distinct delays.}
    \label{fig:6}
\end{figure}

In Figure \ref{fig:5}, the stability recovery of the closed-loop system and the convergence of the plant input to the neighborhood of the optimal value $\theta^{\ast}=(0,1)$ is confirmed. \textcolor{black}{For comparison, Figure \ref{fig:11} illustrates the input trajectory for the case with equal delays, previously discussed in Section \ref{Sec4.2}. Notably, the response in the equal-delay scenario is slightly smoother, likely due to the symmetric compensation applied to both control variables. Nevertheless, in both cases, the Newton-based ESC successfully ensures stabilization around the optimal point, reinforcing the effectiveness of the proposed approach.} Furthermore, Figure \ref{fig:6} depicts that the estimate of $\hat{\theta}(t)$ also converges to the optimum value of the nonlinear map, despite the presence of time delays.

\begin{figure}[htb!]
    \centering    \includegraphics[scale=0.54]{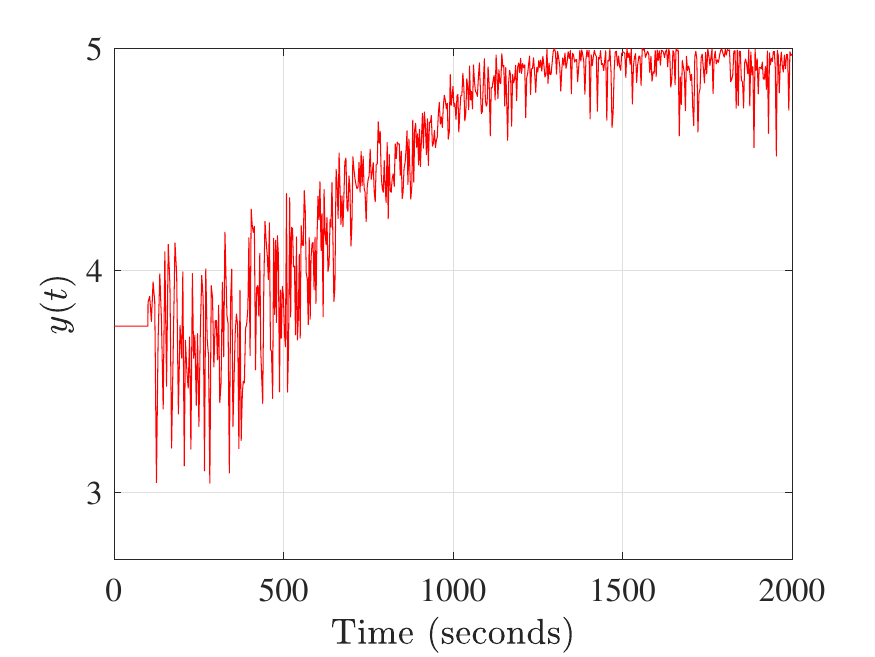}  \caption{Evolution of the system output $y(t)$ under the multivariable Newton-based stochastic ESC and predictor feedback when input channels experience distinct delays.}
    \label{fig:7}
\end{figure}
\begin{figure}[htb!]
    \centering    \includegraphics[scale=0.54]{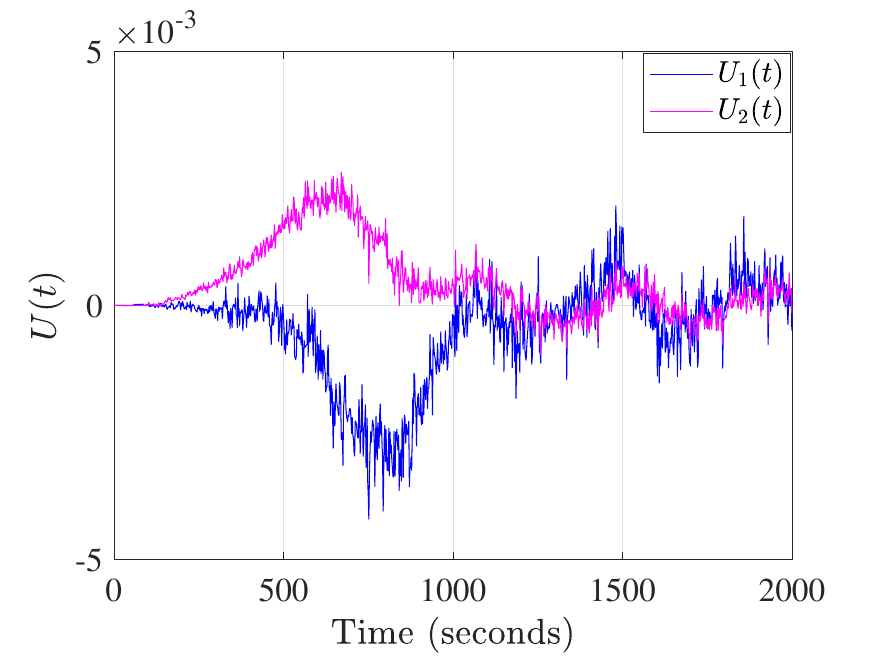}
    \caption{Control signal $U(t)$ for the multivariable Newton-based stochastic ESC and predictor feedback with distinct delays.}
    \label{fig:8}
\end{figure}
\begin{figure}[htb!]
    \centering
    \includegraphics[scale=0.54]{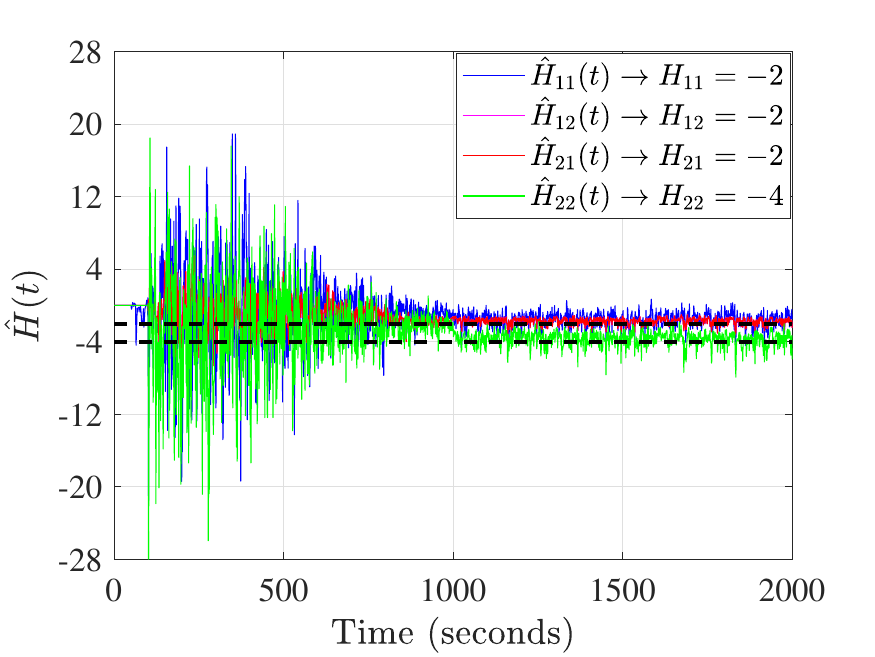}
    \caption{Convergence of the elements of $\hat{H}(t)$ to the Hessian matrix $H$. The first black dashed line corresponds to $H_{11}=H_{12}=H_{21}=-2$, and the second black dashed line corresponds to $H_{22}=-4$.}
    \label{fig:9}
\end{figure}
\newpage
\begin{figure}[htb!]
    \centering    \includegraphics[scale=0.54]{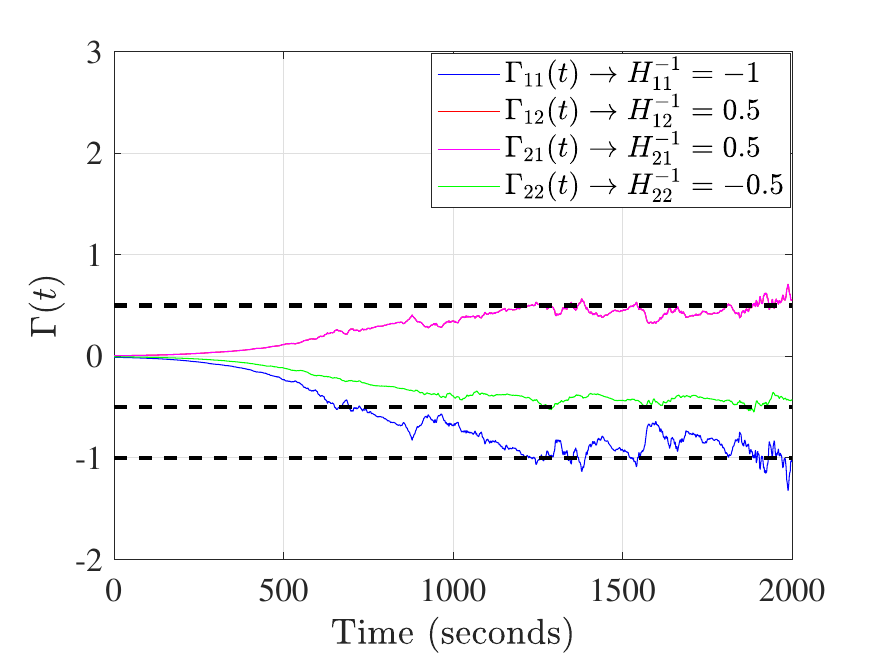}  \caption{Convergence of $\Gamma(t)$ to $H^{-1}$. The first black dashed line corresponds to $H^{-1}_{11}=-1$, the second black dashed line corresponds to $H^{-1}_{12}=H^{-1}_{21}=0.5$, and the third black dashed line corresponds to $H^{-1}_{22}=-0.5$.}
    \label{fig:10}
\end{figure}

Figure \ref{fig:7} indicates that the total delay introduced into the system was effectively compensated by the predictor feedback. As expected, the plant output converged to the desired value of $y^{\ast}=5$.

In Figure \ref{fig:8}, there is a noticeable attenuation of the control signal $U(t)$ as the output $y(t)$ of the map tends towards the extremum point. From Figures \ref{fig:9} and \ref{fig:10}, it is observed that the estimates of the Hessian matrix (\ref{eq:Hessiana}) and its inverse, by means of (\ref{eq:15}), converge to the unknown elements of $H$ and $H^{-1}$, respectively.

\textcolor{black}{
Increasing the adaptation gain $K$ enhances the closed-loop response speed \cite{IJRNC2011,IJACSP2016}, while reducing the constants $a_i$ decreases the amplitude of residual oscillations. As established in Theorem~1, the ultimate residual set for the error $\theta(t)-\theta^*$ is of the order $O(|a| + 1/\omega)$. }




\newpage

\section{Conclusion}\label{sec13}
This paper introduced a novel multivariable Newton-based stochastic ESC scheme, developed for real-time optimization in systems with distinct input delays across multiple actuator channels. By employing predictor-based feedback with Hessian inverse estimation and carefully tuned stochastic excitation signals, the proposed control law effectively compensates for delays, improving the system’s robustness and adaptability in dynamic environments. 

This approach was extended to manage systems with multiple input channels and a single output, effectively handling channel cross-coupling and preserving exponential stability. Through rigorous stability analysis based on backstepping transformations and the averaging theorem in infinite dimensions, the scheme guarantees convergence of the output to a small neighborhood around the extremum, even under arbitrarily long input delays. Numerical results validated the theoretical framework, highlighting both the complexities and the advantages of addressing delay compensation in real-time control.

\textcolor{black}{Our approach can also accommodate known time-varying delays and even unknown time-varying delays that evolve sufficiently slowly. However, the averaging results for functional differential equations \cite{AHale} assume constant delays, which imposes a theoretical limitation. Additionally, our method can be adapted to handle unknown constant delays, as demonstrated in \cite{ref_35_TAC}. Nevertheless, due to the nonlinear parametrization of the delay, this adaptation would only guarantee local convergence --- requiring the initial delay estimate to be sufficiently close to the actual delay. This localized approach offers limited advantages over the inherent robustness of our predictor-feedback strategy to small delay perturbations, as rigorously analyzed in \cite{ref_36_TAC}.}

A fundamental attribute of ESC emphasized here is its design for real-time optimization of unknown objective functions. This capability makes ESC an ideal choice for applications requiring immediate, adaptive responses without prior knowledge of the objective landscape --- an advantage that sets it apart from offline optimization methods focused on known objectives, with or without constraints \citep{Opt1, Opt2, Opt3}.

In light of the rapid advancements in learning-based algorithms \citep{Learn1, Learning2, Learning3} --- particularly for optimization and control --- ESC offers a robust, learning-based, model-free, and data-driven \citep{Model1, Data1, Data2, Data3, Model2} alternative that combines real-time processing with rigorous convergence guarantees. 
Unlike many machine learning and reinforcement learning algorithms, which often lack consistent, provable convergence standards, ESC offers real-time adaptability and reliability essential for critical control applications. This positions ESC as a tool capable of bridging the adaptability of learning-based approaches with the precision and assurance needed in real-time systems.

Overall, the presented work broadens the applicability of ESC by addressing multivariable systems with distinct input delays, paving the way for further exploration in complex, delay-sensitive environments. This framework could be particularly impactful in applications where real-time performance, stability under delay, and minimal reliance on predefined models are crucial --- such as in autonomous systems, robotics, and adaptive industrial processes \cite{AJC2014}.


\textbf{Future Directions:} 

\textcolor{black}{The time delays considered in this work can be arbitrarily large, as they are assumed to be known. However, the proposed predictor-feedback strategy is robust to small delay mismatches, as discussed in \cite{ref_353_HP_Krstic}. Furthermore, robustness to measurement noise is ensured by averaging theory for classes of stochastic and Gaussian noise (with zero mean), as such disturbances are naturally filtered out in the feedback loop.} 

\textcolor{black}{An alternative bounded extremum-seeking framework was proposed in \cite{ref_186_HP_Krstic}, where the control output maintains analytically guaranteed bounds despite operating on unknown, noisy, and time-varying functions.
Additionally, \citet{TNM:2006} analyzed several extremum-seeking schemes, demonstrating that, under appropriate conditions, these approaches achieve extremum seeking from arbitrarily large initial condition domains when controller parameters are properly tuned.}

\textcolor{black}{Another promising research direction involves applying the heuristic Luus-Jaakola algorithm \cite{Luus_Jaakola}, which does not require convexity, differentiability, or even local Lipschitz continuity, making it particularly suited for non-differentiable and non-convex optimization problems.}

\textcolor{black}{From a multi-output (payoff-function) perspective, extremum-seeking problems are closely related to game theory, particularly Nash Equilibrium Seeking (NES). See \cite[Part III]{Oliveira2022} for NES algorithms that account for delays and PDEs. However, the integration of stochastic extremum-seeking into NES under delays remains an open challenge, making it a promising direction for future research.}

\textcolor{black}{Another compelling challenge is the optimal design of perturbation signals $S(t)$ for extremum-seeking, as an alternative to modifying the demodulation signals $M(t)$ and $N(t)$. This problem, commonly called ``trajectory generation'' or ``motion planning,'' involves designing perturbation signals that anticipate their propagation through PDEs, extending beyond transport PDEs or pure delay systems. In particular, this requires solving Cauchy-Kovalevskaya-type motion planning problems for systems governed by parabolic and hyperbolic PDEs \cite[Chapter~12]{livro_azul_Krstic}.}


\textcolor{black}{However, one could explore how Physics-Informed Neural Networks (PINNs) or other machine-learning techniques could automate this process for a broader class of PDEs, eliminating the need for case-by-case analytical derivations. We believe that integrating extremum-seeking control with machine learning offers a promising and innovative approach, enabling the experimental discovery of equivalent solutions through PINNs. This fusion could significantly extend the applicability of extremum-seeking to even more complex dynamical systems while preserving rigorous convergence guarantees.} \\

\backmatter





\bmhead{Acknowledgements}
This work was supported in part by the Brazilian funding agencies CAPES - finance code 001, CNPq, and FAPERJ.

\bibliography{sn-bibliography}

\end{document}